\documentclass[a4paper,12pt]{amsart}

\usepackage{amssymb}
\usepackage{amsmath}
\usepackage{color}
\usepackage{mathrsfs}

 \newtheorem{theorem}{Theorem}
 
 \newtheorem{lemma}[theorem]{Lemma}
 \newtheorem{proposition}[theorem]{Proposition}
 \theoremstyle{definition}
 
 \theoremstyle{remark}
 \newtheorem{remark}[theorem]{Remark}
 \newtheorem{example}[theorem]{Example}

\numberwithin{equation}{section}
\numberwithin{theorem}{section}

\newcommand{\C}{\mathbb C}
\newcommand{\R}{\mathbb R}
\newcommand{\N}{\mathbb N}
\newcommand{\Disc}{\mathbb D}

\newcommand{\D}{\mathcal{D}}
\newcommand{\h}{\mathcal{H}}
\newcommand{\m}{\mathcal{M}}
\newcommand{\M}{\mathcal{M}}
\newcommand{\A}{\mathcal{A}}
\newcommand{\E}{\mathcal{E}}

\newcounter{obs}

\begin{document}

\title[The Ces\`{a}ro space of Dirichlet series]
{The Ces\`{a}ro   space of Dirichlet series \\
and its multiplier algebra}

\author[J.\ Bueno-Contreras]{J.\ Bueno-Contreras}
\address{Departamento de An\'alisis Matem\'atico  \& IMUS,  Universidad de Sevilla,
Calle Tarfia s/n,  Sevilla 41012, Spain}
\email{jjbueno@us.es}

\author[G.\ P.\ Curbera]{G.\ P.\ Curbera}
\address{Facultad de Matem\'aticas \& IMUS,  Universidad de Sevilla,
Calle Tarfia s/n,  Sevilla 41012,  Spain}
\email{curbera@us.es}

\author[O.\ Delgado]{O.\ Delgado}
\address{Departamento de Matem\'atica Aplicada I, E.\ T.\ S.\ de Ingenier\'ia
de Edificaci\'on, Universidad de Sevilla, Avenida de Reina Mercedes
4 A,  Sevilla 41012, Spain}
\email{olvido@us.es}

\thanks{The authors acknowledge the support of MTM2015-65888-C4-1-P,
MINECO (Spain).}

\thanks{This work is part of the PhD thesis of Jos\'e Jorge Bueno Contreras, defended
in the University of Sevilla under the direction of Guillermo P. Curbera
and Olvido Delgado.}


\date{\today}


\begin{abstract}
We consider the space  $\h(ces_p)$ of all Dirichlet series
whose coefficients belong to the Ces\`{a}ro sequence space $ces_p$,
consisting of all complex sequences  whose absolute Ces\`{a}ro means
are in $\ell^p$, for $1<p<\infty$.
It is a Banach space of analytic functions, for which
we study the maximal domain of analyticity and the boundedness of point evaluations.
We identify the algebra of analytic multipliers on $\h(ces_p)$ as the
Wiener algebra of Dirichlet series shifted to the vertical
half-plane $\C_{1/q}:=\{s\in\C:\Re s>1/q\}$, where $1/p+1/q=1$.
\end{abstract}

\maketitle


\section{Introduction}\label{Intro}


Several spaces of Dirichlet series have been studied in recent years.
Hedenmalm, Lindqvist and Seip  introduced  in \cite{hedenmalm-lindqvist-seip} the
Hilbert space of Dirichlet series $\h$, consisting of all Dirichlet series
\begin{equation*}
f(s):=\sum_{n=1}^\infty a_nn^{-s},\quad s\in\C,
\end{equation*}
with square summable coefficients, $(a_n)_{n=1}^\infty\in\ell^2$.
They used it for solving
a problem discussed by Beurling on complete sequences
in the space $L^2(0,1)$. Due to the Cauchy-Schwarz inequality,
each $f\in\h$ defines  an analytic function on the vertical
half-plane $\C_{1/2}:=\{s\in\C:\Re(s)>1/2\}$.
The space $\h$ becomes a Banach space of analytic functions on
$\C_{1/2}$ when endowed with  the norm
\begin{equation*}
\Vert f\Vert_{\h}:=\Vert(a_n)_{n=1}^\infty\Vert_{\ell^2},\quad f\in\h .
\end{equation*}
The Hardy spaces of Dirichlet series
$\h^p$, for $1\le p<\infty$, were introduced  by Bayart in \cite{bayart}.
They are given by the completion of the space of
Dirichlet polynomials $P(s):=\sum_{n=1}^N a_nn^{-s}$ for the norm
$$
\|P\|_{\h^p}:=\bigg(\lim_{T\to\infty} \frac{1}{2T}\int_{-T}^T|P(it)|^pdt\bigg)^{1/p}.
$$
The space $\h$ corresponds to $\h^p$ for $p=2$.
In \cite{maccarthy}, McCarthy considered the
 weighted Hilbert spaces of Dirichlet series
\begin{equation*}
\h_\alpha:=\bigg\{f(s)=
\sum_{n=2}^\infty a_nn^{-s}: \sum_{n=2}^\infty|a_n|^2(\log n)^{\alpha}<\infty\bigg\},
\end{equation*}
for $\alpha<0$, endowed with the norm
\begin{equation*}
\Vert f\Vert_{\h_\alpha}:=\Vert(a_n(\log n)^{\alpha/2})_{n=2}^\infty
\Vert_{\ell^2},\quad f\in\h_\alpha .
\end{equation*}
More recently, Bailleul and Lef\`{e}vre have studied certain
classes of Bergman-type spaces of Dirichlet series,
$\mathscr{A}_\mu^p$ and $\mathscr{B}^p$, for $1\le p<\infty$, \cite{bailleul-lefevre}.
Another type of weighted Hilbert spaces of Dirichlet series
$\mathscr{D}_\alpha$, for $\alpha>0$, has been considered by Bailleul and Brevig
in \cite{bailleul-brevig}.
It should be highlighted
that  the spaces $\h$, $\h^p$, $\h_\alpha$, $\mathscr{A}_\mu^p$, $\mathscr{B}^p$,
$\mathscr{D}_\alpha$ are  all Banach spaces of analytic functions on the vertical
half-plane $\C_{1/2}$.

A deep feature of Dirichlet series is their product.
The pointwise product $f(s)\cdot g(s)$ of
two Dirichlet series $f(s)=\sum_{n=1}^\infty a_nn^{-s}$ and
$g(s)=\sum_{n=1}^\infty b_nn^{-s}$ is,
in the appropriate domain,
the Dirichlet series $h(s)=\sum_{n=1}^\infty c_nn^{-s}$
whose coefficients $c=(c_n)_{n=1}^\infty$ are given by the \emph{Dirichlet convolution}
$c:=a\cdot b$ of the sequences $a=(a_n)_{n=1}^\infty$ and $b=(b_n)_{n=1}^\infty$, that is,
$$
c_n=(a\cdot b)_n:=\sum_{k|n}a_kb_{\frac{n}{k}} ,\quad n\ge1,
$$
where  $k|n$ denotes that $k$ is a divisor of $n$.

Given a space $\mathcal{E}$ of Dirichlet series, a
\textit{multiplier} on $\mathcal{E}$ is  an analytic function $f$ with the property
that $fg\in \mathcal{E}$ for every $g\in\mathcal{E}$.
The \textit{multiplier algebra} of $\mathcal{E}$ is the set of all
multipliers on $\mathcal{E}$; we denote it by $\m(\E)$.
Note that $\m(\E)\subseteq \E$ whenever $\mathbf{1}\in\E$. Neither of the
spaces $\h$, $\h^p$, $\h_\alpha$, $\mathscr{A}_\mu^p$, $\mathscr{B}^p$, $\mathscr{D}_\alpha$
is closed under multiplication.
Hence,  a relevant question is to identify the \textit{multiplier algebra} of these spaces.
Hedenmalm, Lindqvist and Seip identified the multiplier
algebra $\M$ of the Hilbert space of Dirichlet series $\h$  proving that
$$
\M=\h^\infty,
$$
where $\h^\infty$ is the algebra  of bounded analytic functions  on
$\C_0:=\{s\in\C:\Re(s)>0\}$
which can be represented as a Dirichlet series in some vertical half-plane,
\cite[Theorem 3.1]{hedenmalm-lindqvist-seip}.
This identification was a key step in solving  Beurling's question on complete sequences
in $L^2(0,1)$.
It is noticeable that for all the spaces $\h^p$,  $\h_\alpha$,
$\mathscr{A}_\mu^p$, $\mathscr{B}^p$, $\mathscr{D}_\alpha$
the multiplier algebra is also the  algebra  $\h^\infty$;
\cite[Theorem 7]{bayart},  \cite[Theorem 1.11]{maccarthy},
\cite[Theorem 10.1 and Theorem 11.21]{bailleul}, \cite[Theorem 3]{bailleul-brevig}.


In this paper we consider the space $\h(ces_p)$, for $1<p<\infty$,
of all Dirichlet series $f(s)=\sum_{n=1}^\infty a_nn^{-s}$ such that the sequence
of coefficients $(a_n)_{n=1}^\infty$ belongs to the Ces\`{a}ro
sequence space $ces_p$. The space  $ces_p$ consists of all
complex sequences  whose absolute Ces\`{a}ro means are in $\ell^p$,
that is, complex sequences $(a_n)_{n=1}^\infty$ satisfying
\begin{equation*}\label{cesp}
\Vert (a_n)_{n=1}^\infty\Vert_{ces_p}:=\left(\sum_{n=1}^\infty\bigg(\,\frac{1}{n}
\sum_{k=1}^n|a_k|\,\bigg)^p\right)^{\frac{1}{p}}<\infty.
\end{equation*}
It is a Banach space of sequences, that arises
in a natural way  from Hardy's inequality,
\begin{equation}\label{hardy}
\sum_{n=1}^\infty\bigg(\frac1n\sum_{k=1}^n|a_k|\bigg)^p
\le \left(\frac p{p-1}\right)^p\sum_{n=1}^\infty|a_n|^p,
\end{equation}
which establishes   the boundedness on $\ell^p$ of the Ces\`{a}ro averaging operator,
\cite[Theorem 326]{hardy-littlewood-polya}.
It has been throughly studied by G.\ Bennett, \cite{bennett} and Jagers, \cite{jagers},
see also \cite{astashkin-maligranda} and the references therein.

We define, for  $f(s)=\sum_{n=1}^\infty a_nn^{-s}\in\h(ces_p)$,
\begin{equation}\label{normp}
\Vert f\Vert_{\h(ces_p)}:=\Vert (a_n)_{n=1}^\infty\Vert_{ces_p}.
\end{equation}
With this definition, $\h(ces_p)$ is linearly isometric to $ces_p$.
The aim of this paper is to study $\h(ces_p)$ as  a Banach space of analytic functions,
to  find the maximal common domain of definition of its Dirichlet series,
to study the boundedness and the norm of point evaluations, and to
identify the algebra of analytic multipliers on $\h(ces_p)$.
As we will see, the situation will turn out to be rather different
to that of the previously studied spaces of Dirichlet series.


The paper is organized as follows. Section \ref{Preliminaries} contains preliminary facts on Dirichlet series and
spaces of bounded Dirichlet series.

In Section \ref{section-1} we study
$\h(ces_p)$ as  a Banach space of analytic functions.
From being isometrically isomorphic to $ces_p$, it follows that the
sequence of monomials $\{m^{-s}:m\ge1\}$ forms an unconditional, boundedly complete
and shrinking Schauder basis for $\h(ces_p)$; in particular, $\h(ces_p)$ is reflexive.
We show that all functions in $\h(ces_p)$  are analytic
on the vertical half-plane $\C_{1/q}$, where $1/p+1/q=1$ (Theorem \ref{t-3.3}). We also
study the boundedness on $\h(ces_p)$ of  point
evaluations: $f\mapsto f(s_0)$ for $s_0\in\C_{1/q}$, giving  sharp  estimates for their norm and
the precise order of growth when $\Re(s)$ approaches
the critical value $1/q$ (Theorem \ref{t-3.4}).

Section \ref{section-3} is devoted to identifying  the multiplier
algebra $\m(\h({ces_p}))$ of $\h(ces_p)$. A first result shows that
$$
\A^{1/q}\subseteq\m(\h({ces_p}))\varsubsetneq\h^\infty(\C_{1/q}),
$$
where $\A^{1/q}$ is the space of all Dirichlet series $f(s)=\sum_{n=1}^\infty a_nn^{-s}$
satisfying the condition $\sum_{n=1}^\infty |a_n|n^{-1/q}<\infty$, and $\h^\infty(\C_{1/q})$
is the algebra of bounded analytic functions  on   $\C_{1/q}$
which can be represented as a Dirichlet series (Theorem \ref{t-4.5}).
The result shows that the situation concerning the multiplier
algebra of $\h(ces_p)$ is completely different from that of  other spaces
of Dirichlet series studied previously in the literature:
in this case,  the multiplier algebra will
not coincide with an algebra of bounded Dirichlet series.

The fact that the multiplier algebras of $\h^p$, $\h_\alpha$,
$\mathscr{A}_\mu^p$, $\mathscr{B}^p$, $\mathscr{D}_\alpha$
coincide with $\h^\infty$
is in accordance with--actually, it follows from--the
situation of multipliers for Hardy spaces on unit disc $\Disc$ of $\C$. More
precisely, it follows from
the classical result of Schur identifying the multiplier algebra
of the Hardy space $H^2(\mathbb{D})$, of all
Taylor series with square summable coefficients, with the space
$H^\infty(\mathbb{D})$ of bounded analytic function on
$\mathbb{D}$, \cite[X p.226]{schur}.

In the search of a conjecture to pursue, it is relevant to recall
the situation regarding multipliers of the space  $H(\Disc,ces_p)$ of all Taylor series on
$\mathbb{D}$ with coefficient belonging to $ces_p$.
It was proven by Curbera and Ricker
that the multiplier algebra of $H(\Disc,ces_p))$ is not
$H^\infty(\mathbb{D})$ but a rather smaller algebra, namely,
the Wiener algebra of  all absolutely  convergent Taylor series, which  is  the smallest algebra
inside $H(\Disc,ces_p)$ which contains the  polynomials,
\cite[Theorem 3.1]{curbera-ricker-2013}, \cite[Theorem 4.1]{curbera-ricker-2014}.

The main result of this paper is that
\begin{equation*}\label{eq-1.2}
\m(\h({ces_p}))= \A^{1/q},
\end{equation*}
with equality of norms (Theorem \ref{t-4.8}).  We attempt an explanation of this unexpected result.
Hardy's inequality \eqref{hardy} shows that  $\ell^p$ is continuously included in $ces_p$;
in fact, the inclusion is proper. Even more, $ces_p$  contains sequences
with arbitrarily large terms. Indeed, given \textit{any}  sequence $(a_k)_{k=1}^\infty$ of
complex numbers, there exists a subsequence $(e^{m_k})_{k=1}^\infty$ of the canonical vectors
$\{e^m:m\ge1\}$ in $\C^\N$ such that
$\sum_{k=1}^\infty a_k e^{m_k}$ belongs to $ces_p$.
This is an important feature of $ces_p$.
Thus, the space $\h(ces_p)$ contains Dirichlet series whose
coefficients can be arbitrarily large. This feature of $\h(ces_p)$
may explain the multiplier algebra being the smallest possible
algebra which contains the  Dirichlet polynomials.

We end in Section \ref{section-4} with two further results about multipliers on $\h(ces_p)$:
regarding  compact multipliers (Theorem \ref{t-5.1}) and on the Schur point-wise
multipliers from $\h(ces_p)$
to $\A^{1/q}$ (Theorem \ref{t-5.2}).


\section{Preliminaries}\label{Preliminaries}


We collect some general facts on Dirichlet series.
Recall that if a Dirichlet series is convergent (or absolutely
convergent) at a point $s_0\in\C$, then it is convergent  (or absolutely convergent) at any point
$s\in\C$ such that $\Re(s)>\Re(s_0)$.
As a consequence, convergence regions for Dirichlet series
are vertical half-planes $\C_\sigma:=\{s\in\C:\Re(s)>\sigma\}$ for $\sigma\in\R$.
Given a Dirichlet series $f(s)=\sum_{n=1}^\infty a_nn^{-s}$,
its \emph{abscissa of convergence}, denoted by $\sigma_c(f)$, is
 the infimum of all $\sigma\in\R$ such that the series converges on
the vertical half-plane $\C_\sigma$; its \emph{abscissa of absolute convergence}
$\sigma_a(f)$ is
 the infimum of all $\sigma\in\R$ such that the series converges absolutely on
$\C_\sigma$; and its \emph{abscissa of uniform convergence}
$\sigma_u(f)$ is
 the infimum of all $\sigma\in\R$ such that the series converges uniformly on
$\C_\sigma$.
It follows that $-\infty\le \sigma_c(f)\le\sigma_u(f)\le\sigma_a(f)\le+\infty$, and
$\sigma_a(f)-\sigma_c(f)\le1$ if both values are finite.
Bohr proved that $\sigma_a(f)-\sigma_u(f)\le1/2$ (the sharpness
of this inequality is a celebrated theorem of Bohnenblust and Hille,
see \cite[Theorem 5.4.2]{queffelec-q}).
There is a further abscissa associated to a Dirichlet series,
the \emph{abscissa of regularity and boundedness},  $\sigma_b(f)$, which is
 the infimum of all $\sigma\in\R$ such that the function $f(s)=\sum_{n=1}^\infty a_nn^{-s}$
(possibly by analytic continuation from a smaller vertical half-plane)
is analytic and bounded on $\C_\sigma$. Bohr's theorem assures that $\sigma_u(f)=\sigma_b(f)$;
\cite{bohr1913}, see also \cite[Theorem 6.2.3]{queffelec-q}.

We denote by $\D$ the set of all Dirichlet series which
are convergent at some point; this can be equivalently
defined as the set of all Dirichlet series such that the
sequence of its coefficients has, at most,  polynomial growth rate.
Given a Banach space of Dirichlet series $\mathcal{E}\subseteq\D$,
the abscissa of  convergence of
$\mathcal{E}$ is defined by $\sigma_c(\mathcal{E}):=\sup\{\sigma_c(f):f\in\mathcal{E}\}$,
and the abscissa of absolute convergence of $\mathcal{E}$
is $\sigma_a(\mathcal{E}):=\sup\{\sigma_a(f):f\in\mathcal{E}\}$.
In the case when $\sigma_c(\mathcal{E})<\infty$, for every $s_0\in\C_{\sigma_c(\mathcal{E})}$ it is
meaningful to consider  the linear functional $\delta_{s_0}$ on $\mathcal{E}$ given by
point evaluation at $s_0$, that is, $f\in\mathcal{E}\mapsto\delta_{s_0}(f):=f(s_0)\in\C$.

Throughout the paper we will consider
$1<p<\infty$, and $q$ will denote the conjugate exponent of $p$, that is,
$1/p+1/q=1$.

Further notation used in the paper follows.
We denote the set of natural numbers $\{1,2,\dots\}$ by $\N$.
As usual, $\R$ and $\C$ denote the fields of real and complex numbers, respectively.
Given a complex number $s\in\C$, its real part is written as $\Re(s)$ and its imaginary part as
$\Im(s)$. For $\theta\in\R$, the vertical half-plane at the abscissa $\theta$ is denoted by
$\C_\theta:=\{s\in\C:\Re(s)>\theta\}$. The unit disc of the complex plane is $\Disc:=\{z\in\C:|z|<1\}$.
For $\Omega$ a region in $\C$, the space of
all analytic functions on $\Omega$ will be denoted by $H(\Omega)$.
For $k,n\in\N$ we write $k|n$ whenever $k$ is a divisor of $n$.
The integer part of  $x\in\R$,
the largest integer which does not exceed $x$, will be denoted by $\lfloor x\rfloor$.

We write $\zeta$ for the \emph{Riemann zeta-function},
$\zeta(s):=\sum_{n=1}^\infty n^{-s}$, for $\Re(s)>1$.
The constant function  with value one is denoted by $\mathbf{1}$.

Spaces of bounded Dirichlet series play an important role.
We collect some relevant facts on them.
For $r\in\R$, the space $\h^\infty(\C_r)$ consists of
all bounded analytic functions on $\C_r$ which can be represented as
a Dirichlet series in some vertical half-plane, that is,
\begin{equation*}\label{hinf-def}
\h^\infty(\C_r):=\D \cap H^\infty(\C_r).
\end{equation*}

Regarding the abscissa of convergence, we have
$$
\sigma_c(\h^\infty(\C_0))=0 \ \textnormal{ and } \ \sigma_a(\h^\infty(\C_0))=1/2.
$$
The first statement follows from Bohr's theorem, and the second
from a Bohnenblust and Hille's theorem; see, for example,
 \cite[Theorem 1.1.2)]{Balasubramanian-Calado-Queffelec}.
For $\h^\infty(\C_r)$ with  $r\not=0$,  consider the
translation map $\tau_r\colon\D\to\D$ given by
$\tau_r(f)(s):=f(s+r)$, that is,
\begin{equation*}
\tau_r\bigg(\sum_{n=1}^\infty a_nn^{-s}\bigg)=\sum_{n=1}^\infty a_n n^{-(s+r)}
=\sum_{n=1}^\infty (a_nn^{-r}) n^{-s}.
\end{equation*}
The translation $\tau_r$ establishes an isometric
isomorphism between $\h^\infty(\C_r)$ and $\h^\infty(\C_0)$ from which it follows that
\begin{equation}\label{eq-2.1}
\sigma_c(\h^\infty(\C_r))=r \ \textnormal{ and } \ \sigma_a(\h^\infty(\C_r))=r+1/2.
\end{equation}

The space $\h^\infty(\C_r)$ is a linear space which will be endowed with the supremum norm
$$
\Vert f\Vert_{\h^\infty(\C_r)}:=\sup_{s\in\C_r}|f(s)|, \quad f\in\h^\infty(\C_r).
$$
The isometric isomorphism between $\h^\infty(\C_r)$ and $\h^\infty(\C_0)$
allows showing the completeness of
$\h^\infty(\C_r)$ for the supremum norm.
The  result of Hedenmalm, Lindqvist and Seip
states that $\h^\infty(\C_0)$ is isometrically isomorphic to the
multiplier algebra $\M$ of the Hilbert space of Dirichlet series $\h$,
\cite[Theorem 3.1]{hedenmalm-lindqvist-seip}.
Since this last space is complete (for the operator norm) it follows that
$\h^\infty(\C_0)$ is complete for the supremum norm.
Hence,   $\h^\infty(\C_r)$ endowed
with the supremum norm is a Banach space.

For issues related to Dirichlet's series we refer the reader
to \cite{bohr}, \cite{hardy-riesz}, \cite{queffelec-q}, \cite[Ch.IX]{titchmarsh}.


\section{The space of Dirichlet series $\h(ces_p)$}\label{section-1}


The space $\h(ces_p)$, endowed with the norm \eqref{normp},
is a Banach space of Dirichlet series
that inherits its functional properties from the sequence space $ces_p$, as $\h(ces_p)$ and
$ces_p$ are linearly isometric. In particular, we have the following result; see \cite{jagers} and
\cite[Proposition 2.1]{curbera-ricker-2014}.

\begin{proposition}\label{p-3.1}
The following statements hold:
\begin{itemize}
\item[(a)] For every $f(s)=\sum_{n=1}^\infty a_n n^{-s}\in\h(ces_p)$  the Dirichlet polynomials $\sum_{n=1}^N a_n n^{-s}$ converge (as $N\to\infty$) to $f$ in the norm of $\h(ces_p)$. Moreover, from the monotonicity of the norm of $ces_p$,
$$
\|f\|_{\h(ces_p)}=\sup_{N\in\N}\bigg\|\sum_{n=1}^N a_n n^{-s}\bigg\|_{\h(ces_p)}.
$$

\item[(b)] The sequence of monomials $\{m^{-s}:m\ge1\}$ is an unconditional, boundedly complete
and shrinking Schauder basis for $\h(ces_p)$. In particular, $\h(ces_p)$ is reflexive.
\end{itemize}
\end{proposition}

A further approximation for functions in $\h(ces_p)$ is possible. Let  $(p_k)_{k=1}^\infty$
denote the sequence of the prime numbers written in increasing order. For $r\in\N$, let
$$
\N_r:=\bigg\{n\in\N:n=\prod_{i=1}^r p_i^{t_i},\
t_1,\ldots,t_r\ge0\bigg\}.
$$
Consider the map $Q_r$ defined by
$$
f(s)=\sum_{n=1}^\infty a_nn^{-s}\mapsto Q_r(f):=\sum_{n\in\N_r} a_n n^{-s}.
$$
The map $Q_r$ is in fact a projection
$Q_r\colon \h(ces_p)\to \h(ces_p)$.
A remarkable property of the projection $Q_r$ is its multiplicativity, namely,
$Q_r(fg)=Q_r(f)Q_r(g)$,
which holds for any pair of Dirichlet series $f$ and $g$, see \cite[p.157]{queffelec-q}.

Similarly to Proposition \ref{p-3.1}.(a), the following result holds.

\begin{proposition}\label{p-3.2}
For each $f(s)=\sum_{n=1}^\infty a_n n^{-s}$ in $\h(ces_p)$ the Dirichlet series
$\sum_{n\in\N_r}a_n n^{-s}$ converge (as $r\to\infty$)  to $f$ in the norm of $\h(ces_p)$.
Moreover,
$$
\|f\|_{\h(ces_p)}=\sup_{r\in\N}\bigg\|\sum_{n\in\N_r} a_n n^{-s}\bigg\|_{\h(ces_p)}.
$$
\end{proposition}

Let us show that $\h(ces_p)$ is a Banach space of analytic functions. For this,
we determine the abscissa of convergence and the abscissa of absolute convergence of
 $\h(ces_p)$.

\begin{theorem}\label{t-3.3}
Every Dirichlet series $f\in \h(ces_p)$
converges, in fact absolutely, on the vertical half-plane $\C_{1/q}$.
Moreover, the value $1/q$ cannot be improved, that is,
$$
\sigma_c(\h(ces_p))=\sigma_a(\h(ces_p))=1/q.
$$
Consequently, $\h(ces_p)$ is a Banach space of analytic functions on $\C_{1/q}$, which is a
maximal domain.
\end{theorem}

\begin{proof}
Let $f(s)=\sum_{n=1}^\infty a_nn^{-s}\in\h(ces_p)$ with $(a_n)_{n=1}^\infty\in ces_p$. 
Set $r>1/q$. It follows that
\begin{eqnarray*}
\sum_{n=1}^\infty\frac{|a_n|}{n^r}
& \le &
r\sum_{n=1}^\infty|a_n|\sum_{k=n}^\infty\frac{1}{k^{r+1}}
=r\sum_{k=1}^\infty\frac{1}{k^{r+1}}\sum_{n=1}^k|a_n|
\\ & \le &
 r\left(\sum_{k=1}^\infty\frac{1}{k^{rq}}\right)^{1/q}
 \left(\sum_{k=1}^\infty\bigg(\frac{1}{k}\sum_{n=1}^k|a_n|\bigg)^p\right)^{1/p}
\\ & = &
r\zeta(rq)^{1/q}\,\Vert f\Vert_{\h(ces_p)}.
\end{eqnarray*}
Then $\sigma_a(f)\le1/q$ for all $f\in\h(ces_p)$ and so $\sigma_a(\h(ces_p))\le1/q$.

On the other hand, for $r>1/p$ set $f(s):=\sum_{n=1}^\infty1/n^{r+s}$. Note that $f\in \h(ces_p)$ as $(n^{-r})_{n=1}^\infty\in\ell^p\subseteq ces_p$. Since $f(s)=\zeta(r+s)$, it follows that $\sigma_c(f)=1-r$ which tends to $1/q$ as $r\to 1/p$. Thus, $\sigma_c(\h(ces_p))\ge1/q$ and the conclusion follows since $\sigma_c(\h(ces_p))\le\sigma_a(\h(ces_p))$.
\end{proof}


We study next boundedness of  the linear functional $\delta_{s_0}$ on $\h(ces_p)$ given by
evaluation at a point ${s_0}\in \C_{1/q}$:
$$
f\in\h(ces_p)\mapsto \delta_{s_0}(f):=f({s_0})\in\mathbb{C}.
$$
Note, for $s_0=\sigma+it\in\mathbb{C}_{1/q}$ and $f(s)=\sum_{n=1}^\infty a_nn^{-s}\in\h(ces_p)$, that the proof of Theorem \ref{t-3.3}
shows
$$
|\delta_{s_0}(f)|=\bigg|\sum_{n=1}^\infty a_nn^{-s_0}\bigg|\le\sum_{n=1}^\infty|a_n|n^{-\sigma}\le\sigma\zeta(\sigma q)^{1/q}\Vert f\Vert_{\h(ces_p)}.
$$
Thus, $\delta_{s_0}$ belongs to the dual space $\h(ces_p)^*$ of $\h(ces_p)$ with  $\Vert\delta_{s_0}\Vert\le\sigma\zeta(\sigma q)^{1/q}$.

We provide sharp  estimates for the norm $\|\delta_{s_0}\|$, the precise order of growth
when the abscissa approaches from the right the critical value $1/q$, and the asymptotic value when
the abscissa increases to $\infty$.

We require the  dual Banach space of $ces_p$. This space was
isometrically identified
by Jagers, \cite{jagers}. A simpler isomorphic identification
was given by Bennett,  \cite[p.61]{bennett}. Following Bennett
the dual space $ces_p^*$ of $ces_p$
can be identified with the sequence space $d(q)$, for $1/p+1/q=1$,  of all complex
sequences $(b_n)_{n=1}^\infty$
satisfying
\begin{equation*}\label{normbennett}
\Vert (b_n)_{n=1}^\infty\Vert_{d(q)}:=\bigg(\sum_{n=1}^\infty\sup_{k\ge n}|b_k|^q\bigg)^{1/q}<\infty.
\end{equation*}
The action of a sequence $(b_n)_{n=1}^\infty\in d(q)$ as an element of $ces_p^*$ is given by the standard pairing
$$
(a_n)_{n=1}^\infty\in ces_p\mapsto\Big\langle(b_n)_{n=1}^\infty,(a_n)_{n=1}^\infty\Big\rangle:=\sum_{n=1}^\infty a_nb_n.
$$
The equivalence between the norms of $ces_p^*$ and $d(q)$ is given,
for $(b_n)_{n=1}^\infty\in ces_p^*$, by
\begin{equation}\label{eq-3.1}
\frac1q\| (b_n)_{n=1}^\infty\|_{d(q)}
\le \|(b_n)_{n=1}^\infty\|_{ces_p^*} \le (p-1)^{1/p}\| (b_n)_{n=1}^\infty\|_{d(q)}.
\end{equation}
The sequence $(\tilde b_n)_{n=1}^\infty$ defined by
$\tilde b_n:=\sup_{k\ge n}|b_k|$, for $n\ge1$, is
known as  the \textit{least decreasing majorant} of the sequence
$(b_n)_{n=1}^\infty$.
%


\begin{theorem}\label{t-3.4}
For each $s_0=\sigma+it\in \C_{1/q}$ the linear functional
$\delta_{s_0}$ is bounded on $\h(ces_p)$, and
$$
\frac1q\zeta(\sigma q)^{1/q} \le \|\delta_{s_0}\|\le (p-1)^{1/p}\zeta(\sigma q)^{1/q}.
$$
Moreover, there is a value $\sigma_p$, depending
only on $p$, such that $\Vert\delta_{s_0}\Vert=\zeta(p)^{-1/p}$ whenever $\sigma\ge\sigma_p$.
\end{theorem}


\begin{proof}
Let $s_0=\sigma+it\in \C_{1/q}$. For $f(s)=\sum_{n=1}^\infty a_nn^{-s}\in\h(ces_p)$, since $(n^{-s_0})_{n=1}^\infty\in d(q)$, we can write
$$
\delta_{s_0}(f)=f(s_0)=\sum_{n=1}^\infty a_nn^{-s_0}=
\Big\langle\big(n^{-s_0}\big)_{n=1}^\infty,(a_n)_{n=1}^\infty\Big\rangle.
$$
Thus, $\delta_{s_0}$ acting on $\h(ces_p)$ can be identified with the sequence
$(n^{-s_0})_{n=1}^\infty$ acting
on $ces_p$. Since $\h(ces_p)$ and $ces_p$ are
isometric, we have that the norms of $\delta_{s_0}$  as an element of
$\h(ces_p)^*$ and of
$(n^{-s_0})_{n=1}^\infty$ as an element of $ces_p^*$ are equal.
Using Bennett's identification of $ces_p^*$ as
the space $d(q)$, from \eqref{eq-3.1}, it follows that
$$
\frac1q\|(n^{-s_0})_{n=1}^\infty\|_{d(q)}
\le \|\delta_{s_0}\|\le (p-1)^{1/p}\|(n^{-s_0})_{n=1}^\infty\|_{d(q)}.
$$
Note that for sequences
$(b_n)_{n=1}^\infty$ such that the sequence $(|b_n|)_{n=1}^\infty$ is decreasing, we have that
$(b_n)_{n=1}^\infty\in d(q)$ if and only if
$(b_n)_{n=1}^\infty\in \ell^q$, and in this case the norms coincide.
Consequently,
$$
\|(n^{-s_0})_{n=1}^\infty\|_{d(q)}=\|(n^{-s_0})_{n=1}^\infty\|_{\ell^q}=
\bigg(\sum_{n=1}^\infty\frac1{n^{\sigma q}}\bigg)^{1/q}=\zeta(\sigma q)^{1/q}.
$$

In order to prove that $\|\delta_{s_0}\|$ becomes constant when $\sigma=\Re(s_0)$ is
sufficiently large (only depending on $p$)  we require the isometric
identification of $ces_p^*$  given by Jagers, \cite{jagers}.
Namely, for $b=(b_n)_{n=1}^\infty\in ces_p^*$ we have
\begin{align}\label{eq-3.2}
\| (b_n)_{n=1}^\infty \|_{ces_p^*} &= \nonumber \\
& \left(\sum_{n\in D(b)}
\bigg(\frac{|b_{m(n)}|-|b_{m(n+1)}|}{B_{m(n)} -B_{m(n+1)} }
\bigg)^q\big(B_{m(n)}-B_{m(n+1)}\big) \right)^{1/q},
\end{align}
where
\begin{align*}
&B_k:=\sum_{j=k}^\infty1/j^p, \quad k\ge1;
\\
&m(1):=\max\Big\{k\in\mathbb{N}\cup\{\infty\}: |b_k|=\max_{j\ge1}|b_j|\Big\},
\end{align*}
and, for $n\ge1$,
\begin{align*}
m(n+1)  := \max\bigg\{ & k\in\N\cup\{\infty\}\ :\ k>m(n),\
\\ & \
\frac{|b_{m(n)}|-|b_k|}{B_{m(n)} -B_k } =\min_{m(n)<j\le\infty}\
\frac{|b_{m(n)}|-|b_j|}{B_{m(n)} -B_j } \bigg\},
\end{align*}
provided $m(n)$ is defined and finite, else $m(n+1)$ is not defined;
and $D(b)$ is the set of all $k\ge1$ such that $m(k)$ is defined and finite.
It is understood $b_\infty=B_\infty=0$.

Note that if $(|b_n|)_{n=1}^\infty$ is strictly decreasing then $m(1)=1$.
Moreover, $m(2)=\infty$ if the condition
\begin{equation}\label{eq-3.3}
\frac{|b_1|-|b_n|}{B_1-B_n}\ge\frac{|b_1|-|b_\infty|}{B_1-B_\infty}=\frac{|b_1|}{B_1}
\end{equation}
is satisfied for all $n\ge2$. In this case $D(b)=\{1\}$ and
so $\| (b_n)_{n=1}^\infty \|_{ces_p^*} =|b_1|\zeta(p)^{-1/p}$.

For $b=(n^{-s_0})_{n=1}^\infty\in ces_p^*$,
we claim that condition \eqref{eq-3.3} holds provided that
$$
\sigma\ge\sigma_p:=p-1+\frac{\log(p-1)+\log\zeta(p)}{\log2}.
$$

Write \eqref{eq-3.3} for this particular sequence:
$$
\frac{1-\frac1{n^\sigma}}{\sum_{j=1}^{n-1}\frac1{j^p}}\ge\frac1{\sum_{j=1}^\infty\frac1{j^p}},
$$
which is equivalent to
$$
\sum_{j=n}^\infty\frac1{j^p}\ge\frac1{n^\sigma}\zeta(p) .
$$
Since
$$
\sum_{j=n}^\infty\frac1{j^p}\ge\frac1{p-1}\cdot\frac1{n^{p-1}},
$$
it suffices to prove that
$$
\frac1{p-1}\cdot\frac1{n^{p-1}}\ge\frac1{n^\sigma}\zeta(p)
$$
holds for all $n\ge2$. We rewrite this condition
as
$$
n^{\sigma-p+1}\ge(p-1)\zeta(p).
$$
It is clear that for the above
inequality to hold, necessarily we must have $\sigma\ge p-1$. In this case,
the sequence $(n^{\sigma-p+1})_{n=1}^\infty$ is increasing. Thus,
it suffices to check the above inequality for $n=2$:
$$
2^{\sigma-p+1}\ge(p-1)\zeta(p),
$$
that is,
$$
\sigma\ge p-1+\frac{\log(p-1)+\log\zeta(p)}{\log2}=\sigma_p.
$$
Therefore, for $b=(n^{-s_0})_{n=1}^\infty$ with $s_0\in\C_{\sigma_p}$, we have that $m(2)=\infty$ and
so $D(b)=\{1\}$. Hence the sum in \eqref{eq-3.2}
has only one term and
$$
\|\delta_{s_0}\|=\| (n^{-s_0})_{n=1}^\infty \|_{ces_p^*}=
\zeta(p)^{-1/p}.
$$
\end{proof}


\begin{remark}\label{r-3.5}
From the proof of Theorem \ref{t-3.3} and Theorem \ref{t-3.4} actually we have,
for $s_0=\sigma+it\in\C_{1/q}$, that
$$
\|\delta_{s_0}\|\le\min\{\sigma,(p-1)^{1/p}\}\zeta(\sigma q)^{1/q}.
$$
Since $1/q<(p-1)^{1/p}$, as the function $x\mapsto x^x$ is increasing on $(1,\infty)$,
we have that
$$
\min\{\sigma,(p-1)^{1/p}\}=\left\{
\begin{matrix}
	\sigma & \mbox{for} &1/q<\sigma\le (p-1)^{1/p}, \vspace{2mm} \\
	(p-1)^{1/p} & \mbox{for} & \sigma> (p-1)^{1/p}.     \\
\end{matrix}
\right.
$$
\end{remark}


The bounds on the norm of point evaluations in Theorem \ref{t-3.4}
and Remark \ref{r-3.5} can be sharpened for $\h(ces_2)$.
\begin{proposition}\label{p-3.6}
Let $1/2<\Re (s_0)=\sigma\le1$ and $\delta_{s_0}\colon \h(ces_2)\to\C$
be the corresponding
point evaluation functional. Then its norm can be written as
$$
\|\delta_{s_0}\|=\left(\sum_{n=1}^\infty
n^2\bigg(\frac1{n^\sigma}-\frac1{(n+1)^\sigma}\bigg)^2\right)^{1/2},
$$
and the following bounds hold
$$
(2^\sigma-1)\sqrt{\zeta(2\sigma)-1}\le\|\delta_{s_0}\|\le\sigma\sqrt{\zeta(2\sigma)-1}.
$$
\end{proposition}

\begin{proof}
We use the isometric identification of $ces_p^*$
by Jagers for $p=2$.

Let $b=(n^{-s_0})_{n=1}^\infty$. We will prove
that in this case, and for every $m\in\N$, the sequence
\begin{equation}\label{eq-3.4}
\left(\frac{|b_m|-|b_n|}{B_m-B_n }\right)_{n=m+1}^\infty
\end{equation}
is strictly increasing. This condition is precisely
\begin{equation*}\label{bnn}
\frac{ \displaystyle\frac1{m^\sigma} - \frac1{n^\sigma} }
{\displaystyle\sum_{k=m}^\infty\frac1{k^2} - \sum_{k=n}^\infty\frac1{k^2}}
<
\frac{ \displaystyle\frac1{m^\sigma} - \frac1{(n+1)^\sigma} }
{\displaystyle\sum_{k=m}^\infty\frac1{k^2} - \sum_{k=n+1}^\infty\frac1{k^2}},
\end{equation*}
which is equivalent to
\begin{equation}\label{eq-3.5}
\frac{\displaystyle \frac1{m^\sigma}-\frac1{n^\sigma}}
{n^2\bigg(\displaystyle \frac1{n^\sigma} - \frac1{(n+1)^\sigma}\bigg)}
<\sum_{k=m}^{n-1}\frac1{k^2}.
\end{equation}
By applying the mean value theorem to the function $f(x)=x^\sigma$ 
on $(m,n)$ and $(n,n+1)$ we obtain, for $1/2<\sigma\le1$, that
\begin{align}\label{eq-3.6}
\frac{\displaystyle \frac1{m^\sigma}-\frac1{n^\sigma}}
{n^2\bigg(\displaystyle \frac1{n^\sigma} - \frac1{(n+1)^\sigma}\bigg)}
=\frac{(n+1)^\sigma(n^\sigma-m^\sigma)}{n^2m^\sigma\big((n+1)^\sigma-n^\sigma\big)}
\le
\frac{n+1}{n^2m}(n-m).
\end{align}

In order to bound the right-hand side of \eqref{eq-3.5}
we use the following inequality:
\begin{equation*}
\int_m^n\frac{dx}{x^2}+\frac12\bigg(\frac1{m^2}-\frac1{n^2}\bigg)
\le
\sum_{k=m}^{n-1}\frac{1}{k^2},
\end{equation*}
see for instance \cite[p.54]{davis-rabinowitz}.
Since
$$
\int_m^n\frac{dx}{x^2}+\frac12\bigg(\frac1{m^2}-\frac1{n^2}\bigg)
=\bigg(\frac1m-\frac1n\bigg)\bigg(1+\frac1{2m}+\frac1{2n}\bigg),
$$
we have
\begin{equation}\label{eq-3.7}
\bigg(\frac{n-m}{mn}\bigg)\bigg(1+\frac1{2m}+\frac1{2n}\bigg)\le \sum_{k=m}^{n-1}\frac{1}{k^2}.
\end{equation}
Then, \eqref{eq-3.6} and \eqref{eq-3.7} reduce the validity of \eqref{eq-3.5} to
$$
\frac{n+1}{n^2m}(n-m)<\bigg(\frac{n-m}{mn}\bigg)\bigg(1+\frac1{2m}+\frac1{2n}\bigg),
$$
which is true since $m<n$. Thus, \eqref{eq-3.5} holds and so, for every $m\in\N$, the sequence
\eqref{eq-3.4} is strictly increasing.

Hence, for each $n\in\N$ we have that $m(n)=n$.
This implies
that $D(b)=\N$ for $b=(n^{-s_0})_{n=1}^\infty$ and so
\begin{equation*}
\|\delta_{s_0}\|=\|(n^{-s_0})_{n=1}^\infty\|_{ces_2^*}=\left(\sum_{n=1}^\infty n^2\bigg(\frac1{n^\sigma}-\frac1{(n+1)^\sigma}\bigg)^2\right)^{1/2}.
\end{equation*}
Since
$$
n\bigg(\frac1{n^\sigma}-\frac1{(n+1)^\sigma}\bigg)=\frac{1}{(n+1)^\sigma}\,g(n^{-1})
$$
where $g(x)={\displaystyle \frac{(1+x)^\sigma-1}{x}}$ decreases in $(0,\infty)$, we have
$$
\frac{2^\sigma-1}{(n+1)^\sigma}\le
n\bigg(\frac1{n^\sigma}-\frac1{(n+1)^\sigma}\bigg)\le\frac\sigma{(n+1)^\sigma}
$$
and so the bounds for $\Vert\delta_{s_0}\Vert$ follow.
\end{proof}


In the case $p=2$, there are two equivalent expressions for
the norm in $ces_2$ (and so for the norm in $\h(ces_2)$) which are of independent interest.

\begin{proposition}\label{p-3.7}
Let $a=(a_n)_{n=1}^\infty\in ces_2$. Define the functionals
\begin{align*}
M(a)&:=\bigg(\sum_{i,j=1}^\infty\frac{|a_i||a_j|}{\max\{i,j\}}\bigg)^{1/2},
\\
N(a)&:=\bigg(\sum_{n=1}^\infty\frac{|a_n|}{n}
\sum_{k=1}^n|a_k|\bigg)^{1/2} .
\end{align*}
Then
$$
N(a)\le M(a) \le \|a\|_{ces_2}\le \sqrt2 M(a)\le 2 N(a).
$$
\end{proposition}

\begin{proof}
Rearranging the sums in the norm of $a$ we obtain that
\begin{align*}
\| a\|_{ces_2}^2 &=
\sum_{n=1}^\infty\bigg(\frac1n\sum_{k=1}^n|a_k|\bigg)^2
\\   &= \sum_{n=1}^\infty\frac1{n^2}\bigg(\sum_{1\le i,j\le n}|a_i||a_j|\bigg)
\\   &=
\sum_{i,j=1}^\infty|a_i||a_j|\bigg(\sum_{n\ge i,j}\frac1{n^2}\bigg).
\end{align*}
Since $1/n\le\sum_{k=n}^\infty k^{-2}\le2/n$ for every $n\ge1$, it follows that
$$
\frac{1}{\max\{i,j\}} \le \sum_{n\ge i,j}\frac1{n^2}
\le \frac{2}{\max\{i,j\}}.
$$
Hence,  we deduce that $M(a)\le\| a\|_{ces_2}\le\sqrt2 M(a)$.

On the other hand,
\begin{align*}
M(a)^2 &= \sum_{i,j=1}^\infty\frac{|a_i||a_j|}{\max\{i,j\}}
= \sum_{n=1}^\infty\frac1n\bigg(\sum_{\max\{i,j\} =n}|a_i||a_j|\bigg)
\\ & =
\sum_{n=1}^\infty\frac{|a_n|}n\bigg(|a_n|+2\sum_{k=1}^{n-1}|a_k|\bigg)
\\ &\le
2\sum_{n=1}^\infty\frac{|a_n|}{n}\sum_{k=1}^{n}|a_k|
= 2N(a)^2.
\end{align*}
In a similar way
$$
M(a)^2
\ge
\sum_{n=1}^\infty\frac{|a_n|}{n}\sum_{k=1}^{n}|a_k|
=N(a)^2.
$$
Consequently, $N(a)\le M(a)\le \sqrt2 N(a)$.
\end{proof}


\section{The multiplier algebra of $\h(ces_p)$}\label{section-3}


Given a Banach space of Dirichlet series $\E\subseteq\D$ with
convergence abscissa $\sigma_c(\E)$, a \textit{multiplier} on $\E$ is
an analytic function $f$ on $\C_{\sigma_c(\E)}$ with the property that $fg\in \E$ for every $g\in\E$.
The \textit{multiplier algebra} of $\E$ is the space of all multipliers on $\E$,
which will be denoted by $\M(\E)$.
Standard arguments give the following facts on $\M(\E)$.

\begin{proposition}\label{p-4.1}
Let $\E\subseteq\D$ be a Banach space of Dirichlet series.
Suppose that there exists $\sigma\ge\sigma_c(\E)$
such that the point evaluation functional $\delta_{s_0}$ is continuous on
$\E$ for every $s_0\in\mathbb{C}_\sigma$.
Then the following holds:
\begin{itemize}
\item[(a)] For every $f\in\M(\E)$, the operator
$M_{f}\colon\E\to\E$,
given by $M_{f}(g):=fg$ for all $g\in\E$,
is linear and bounded.
\item[(b)] If the constant function $\mathbf{1}\in\E$, then
$\M(\E)\subseteq \E$ and for every $f\in\M(\E)$ it follows that
$\Vert f\Vert_{\E}\le
\Vert \mathbf{1}\Vert_{\E}\Vert f\Vert_{\M(\E)}$,
where $\| f\|_{\M(\E)}$ denotes the operator norm of $M_f$.
Moreover, in this case $\M(\E)$ is a closed
subspace of the space $B(\E)$ of
all bounded linear operators of $\E$ into itself, and,
endowing $\M(\E)$ with the norm $\| \cdot\|_{\M(\E)}$,
the inclusion $\M(\E)\subseteq \E$ is continuous with embedding
constant equal to $\|\mathbf{1}\|_{\E}$.
\end{itemize}
\end{proposition}


The next proposition shows that, under minimal conditions which guarantee a good
behavior of $\M(\E)$, every multiplier on $\E$ is  a bounded analytic function on the appropriate domain.


\begin{proposition}\label{p-4.2}
Let $\E\subseteq\D$ be a Banach space of Dirichlet series
satisfying the condition of Proposition \ref{p-4.1} for some $\sigma\ge\sigma_c(\E)$ and such that
$\mathbf{1}\in\E$.
Then,
$$
\M(\E)\subseteq\h^\infty(\C_\sigma),
$$
where the inclusion is continuous with continuity constant equal to one.
\end{proposition}

\begin{proof}
Let $f\in\M(\E)$. By Proposition \ref{p-4.1}, we have that
$\M(\E)\subseteq\E$ and so
$f^2=ff\in\E$ with
$$
\|f^2\|_{\E}
\le
\|f\|_{\E} \|f\|_{\M(\E)}
\le
\| \mathbf{1} \|_{\E}\|f\|_{\M(\E)}^2.
$$
Iterating the above procedure, we obtain, for every $n\ge1$,  that $f^n\in \E$ and
$$
\| f^n\|_{\E}\le\| \mathbf{1} \|_{\E}\| f\|_{\M(\E)}^n.
$$
For each $s_0\in\C_{\sigma}$, by hypothesis, the point evaluation functional
$\delta_{s_0}$ is bounded on $\E$.
Then
$$
|f^n(s_0)|=|\delta_{s_0}(f^n)|
\le
\|\delta_{s_0}\|\cdot\| f^n\|_{\E}
\le
\|\delta_{s_0}\|\cdot\| \mathbf{1} \|_{\E}\| f\|_{\M(\E)}^n.
$$
Since $|f^n(s_0)|=|f(s_0)|^n$, it follows that
$$
|f(s_0)|\le\big(\|\delta_{s_0}\|\cdot\| \mathbf{1} \|_{\E}\big)^{1/n}\| f\|_{\M(\E)}.
$$
Making $n\to\infty$ we have that $|f(s_0)| \le \| f\|_{\M(\E)}$.
Hence, $f\in \h^\infty(\C_{\sigma})$ and $\Vert f\Vert_{\h^\infty(\C_{\sigma})}\le\| f\|_{\M(\E)}$.
\end{proof}


Furthermore, if the monomials $n^{-s}$, for $n\ge1$, are
multipliers on $\E$, a certain natural weighted $\ell^1$-space of
Dirichlet series is included in $\M(\E)$.

\begin{proposition}\label{p-4.3}
Let $\E\subseteq\D$ be a Banach space of Dirichlet series satisfying the
condition of Proposition \ref{p-4.1} for some $\sigma\ge\sigma_c(\E)$ and such that
$\mathbf{1}\in\E$. Suppose that  $\{n^{-s}: n\ge1\}\subset\M(\E)$ and
denote $\mu_n:=\Vert n^{-s}\Vert_{\M(\E)}$ for $n\ge1$. Then
$$
\mathcal{A}((\mu_n)_{n=1}^\infty):=\bigg\{f(s)
=\sum_{n=1}^\infty a_nn^{-s}:\ \sum_{n=1}^\infty|a_n|\mu_n<\infty\bigg\}\subseteq\M(\E)
$$
and $\Vert f\Vert_{\M(\E)}\le\sum_{n=1}^\infty|a_n|\mu_n$
for all $f\in\mathcal{A}((\mu_n)_{n=1}^\infty)$.
\end{proposition}

\begin{proof}
Let $f(s)=\sum_{n=1}^\infty a_nn^{-s}\in\A((\mu_n)_{n=1}^\infty)$.
The series $\sum_{n=1}^\infty a_nn^{-s}$ is absolutely convergent in $\M(\E)$, as
$$
\sum_{n=1}^\infty\Vert a_nn^{-s}\Vert_{\M(\E)}=\sum_{n=1}^\infty|a_n|\mu_n<\infty,
$$
and so it converges in norm to some $h\in\M(\E)$.
Since, $\M(\E)\subseteq\E$ continuously and so norm convergence in
$\M(\E)$ implies pointwise convergence on $\C_\sigma$, it follows
that $f=h\in \M(\E)$.
From the equality above it follows
that $\|f\|_{\M(\E)}\le\sum_{n=1}^\infty|a_n|\mu_n$.
\end{proof}


\begin{remark}\label{r-4.4}
The particular spaces $\mathcal{A}((\mu_n)_{n=1}^\infty)$ above obtained
for $r\in\R$ and $\mu_n:=n^{-r}$ for all $n\ge1$, are denoted by
\begin{equation*}\label{ar-def}
\A^r:=\bigg\{f(s)=\sum_{n=1}^\infty a_nn^{-s}:\ \sum_{n=1}^\infty |a_n|n^{-r}<\infty\bigg\}.
\end{equation*}
They are Banach spaces when endowed with the norm
$\Vert f\Vert_{\A^r}:=\sum_{n=1}^\infty |a_n|n^{-r}$.
When $r=0$, the corresponding space is the well known \emph{Wiener-Dirichlet algebra}
$\A^+$, see \cite{bayart-finet-li-Queffelec}.
Direct computation shows that $\sigma_c(\A^r)=\sigma_a(\A^r)=r$.
Regarding the point evaluations on $\mathcal{A}^r$,
we have that $\Vert\delta_{s_0}\Vert=1$, for every $s_0\in\mathbb{C}_r$.
With respect to the multipliers, by Proposition \ref{p-4.1},
$\m(\A^r)\subseteq \A^r$ continuously
with embedding constant equal to one.
In fact, both spaces coincide with equality of norms. To see this,
we check that monomials are multiplier on $\A^r$.  For $m\in\N$,
consider $m^{-s}$ and let $g(s)=\sum_{n=1}^\infty a_nn^{-s}\in\A^r$. 
Noting that $m^{-s}g(s)=\sum_{n=1}^\infty c_nn^{-s}$ with 
$c_n=a_{\frac{n}{m}}$ if $m|n$ and $c_n=0$ in other case, it follows
$$
\sum_{n=1}^\infty|c_n|n^{-r}=\sum_{{n=1 \atop m|n}}^\infty|a_{\frac{n}{m}}|n^{-r}
=\sum_{k=1}^\infty|a_k|(km)^{-r}
=m^{-r}\Vert g\Vert_{\A^r},
$$
and so $m^{-s}g\in\A^r$ with $\Vert m^{-s}g\Vert_{\A^r}=m^{-r}\Vert g\Vert_{\A^r}$.
Hence, $m^{-s}\in\M(\A^r)$ and $\Vert m^{-s}\Vert_{\M(\A^r)}\le m^{-r}$.
Actually $\Vert m^{-s}\Vert_{\M(\A^r)}=m^{-r}$, as $\mathbf{1}\in\A^r$
with norm one and $\Vert m^{-s}\Vert_{\A^r}=m^{-r}$. Therefore, from
Proposition \ref{p-4.3}, the inclusion $\A^r\subseteq\M(\A^r)$ holds continuously
with embedding constant equal to one.
\end{remark}


We now consider the space $\m(\h(ces_p))$ of all multipliers  on $\h(ces_p)$.
We have seen that $\sigma_c(\h(ces_p))=1/q$ and that, for each $s_0\in\mathbb{C}_{1/q}$,
the point evaluation functional $\delta_{s_0}$ is continuous on $\h(ces_p)$,
see Theorem \ref{t-3.3} and Theorem \ref{t-3.4}.
With these conditions, Proposition \ref{p-4.1} implies that
every multiplier $f$ on $\h(ces_p)$ defines a bounded multiplication operator $M_f$
from $\h(ces_p)$ into itself:
$$
g\in\h(ces_p)\mapsto M_f(g)=fg\in\h(ces_p).
$$
Moreover, since the constant function $\mathbf{1}\in\h({ces_p})$ and
$\| \mathbf{1}\|_{\h({ces_p})}={\zeta(p)}^{1/p}$, it also follows that
\begin{equation*}
\M(\h(ces_p))\subseteq \h(ces_p)
\end{equation*}
continuously with embedding constant ${\zeta(p)}^{1/p}$, that is,
\begin{equation*}
\Vert f\Vert_{\h({ces_p})}\le{\zeta(p)}^{1/p}\cdot\Vert f\Vert_{\m(\h({ces_p}))},
\quad f\in\m(\h({ces_p})).
\end{equation*}

So, a multiplier $f$ on $\h(ces_p)$ is actually a Dirichlet series $f(s)=\sum_{n=1}^\infty a_nn^{-s}$ belonging to $\h(ces_p)$ and the action of the multiplication operator $M_f$ on $g(s)=\sum_{n=1}^\infty b_nn^{-s}\in\h(ces_p)$ is given by
$$
M_f(g)(s)=f(s)g(s)=\sum_{n=1}^\infty \bigg(\sum_{k|n}a_kb_{\frac{n}{k}}\bigg)n^{-s}.
$$
The boundedness of the operator  $M_f$ corresponds to
the existence of some constant $M>0$ such that
$$
\left(\sum_{n=1}^\infty\bigg(\frac1n\sum_{k=1}^n\bigg|
\sum_{j|k}a_jb_{\frac{k}{j}}\bigg|\bigg)^p \right)^{1/p}
\le
M\cdot
\left(\sum_{n=1}^\infty\bigg(\frac1n\sum_{k=1}^n|b_k|\bigg)^p\right)^{1/p}
$$
for all $g(s)=\sum_{n=1}^\infty b_nn^{-s}\in\h(ces_p)$.
The least of such constants $M$
is the operator norm $\|M_f\|$ of $M_f$ as a bounded  operator from $\h(ces_p)$ into itself,
which we denote by $\| f\|_{\M(\h({ces_p}))}$.


\begin{theorem}\label{t-4.5}
The inclusions
$$
\A^{1/q}\subseteq\m(\h({ces_p}))\varsubsetneq\h^\infty(\C_{1/q})
$$
holds continuously with inclusion constants equal to one.
\end{theorem}

\begin{proof}
The first inclusion follows from Proposition \ref{p-4.3}
if we show, for $m\ge1$, that the monomial $m^{-s}$ is a multiplier on $\h(ces_p)$,
and
\begin{equation*}
\big\Vert m^{-s}\big\Vert_{\m(\h({ces_p}))}=m^{-1/q}.
\end{equation*}

Let $g(s)=\sum_{n=1}^\infty b_nn^{-s}\in\h(ces_p)$. The coefficients of the
Dirichlet series $m^{-s}g$ are given by
$$
(m^{-s}\cdot g)_k=
\left\{\begin{array}{ll}
b_i & \textnormal{ if } k=im \textnormal{ for some } i\ge1 \\ 0 & \textnormal{ in other case}
\end{array}\right..
$$

We estimate:
\begin{eqnarray*}
\|m^{-s}g \|_{\h(ces_p)}^p
&=&
\sum_{n=1}^\infty\bigg(\frac{1}{n}\sum_{k=1}^n|(m^{-s}\cdot g)_k|\bigg)^p
=\sum_{n=m}^\infty\bigg(\frac{1}{n}\sum_{i=1}^{\lfloor n/m\rfloor}|b_i|\bigg)^p
\\ & = &
\sum_{j=1}^\infty\sum_{n=jm}^{(j+1)m-1}\bigg(\frac{1}{n}\sum_{i=1}^j|b_i|\bigg)^p
=\sum_{j=1}^\infty\bigg(\sum_{i=1}^j|b_i|\bigg)^p\,\sum_{n=jm}^{(j+1)m-1}\frac{1}{n^p}
\\ & \le &
\sum_{j=1}^\infty\bigg(\sum_{i=1}^j|b_i|\bigg)^p\,\frac{m}{(jm)^p}
=
\frac{1}{m^{p-1}}\,\Vert g\Vert_{\h({ces_p})}^p .
\end{eqnarray*}
Then $\Vert m^{-s}g\Vert_{\h({ces_p})}\le m^{-1/q}\Vert g\Vert_{\h({ces_p})}$.
Thus, $m^{-s}$ is a multiplier on $\h({ces_p})$ and
$ \Vert m^{-s}\Vert_{\m(\h({ces_p}))}\le m^{-1/q}$.

On the other hand, for $g(s)=j^{-s}$, with $j\ge2$, we have
$$
\Vert m^{-s} j^{-s}\Vert_{\h({ces_p})}^p
=\sum_{n= jm}^\infty\frac{1}{n^p}\ge\frac{1}{p-1} \frac{1}{(jm)^{p-1}}
$$
and
$$
\Vert j^{-s}\Vert_{\h({ces_p})}^p
=\sum_{n=j}^\infty\frac{1}{n^p}\le\frac{1}{p-1}\frac{1}{(j-1)^{p-1}}.
$$
Hence,
$$
\Vert m^{-s}\Vert_{\m(\h({ces_p}))}
\ge
\frac{\Vert m^{-s}\cdot j^{-s}\Vert_{\h({ces_p})}}{\Vert j^{-s}\Vert_{\h({ces_p})}}
\ge
\frac{(j-1)^{1/q}}{(jm)^{1/q}}.
$$
Making $j\to\infty$, we arrive at $\Vert m^{-s}\Vert_{\m(\h({ces_p}))}\ge m^{-1/q}$.

The second inclusion follows from  Proposition \ref{p-4.2}.

It only remains to prove that $\M(\h({ces_p}))\not=\h^\infty(\C_{1/q})$.
For this we calculate the abscissa
of convergence and absolute convergence of $\M(\h({ces_p}))$.
From Theorem \ref{t-3.3}, Remark \ref{r-4.4} and the inclusions
$$
\A^{1/q}\subseteq\M(\h({ces_p}))\subseteq\h({ces_p}),
$$
it follows that
\begin{align*}
1/q=\sigma_c(\A^{1/q}) &\le\sigma_c\big(\M(\h({ces_p}))\big)
\\
&\le
\sigma_a\big(\M(\h({ces_p}))\big)\le\sigma_a(\h({ces_p}))=1/q.
\end{align*}
Then $\sigma_c\big(\M(\h({ces_p}))\big)=\sigma_a\big(\M(\h({ces_p}))\big)=1/q$. Thus, $\M(\h(ces_p))\not=\h^\infty(\mathbb{C}_{1/q})$ as
$\sigma_a(\h^\infty(\C_{1/q}))=1/q+1/2$ by \eqref{eq-2.1}.
\end{proof}


\bigskip


Theorem \ref{t-4.5}  already shows that the situation concerning the multiplier
algebra of $\h(ces_p)$ is certainly different from that of  other spaces
of Dirichlet series studied previously in the literature.
In this case,  the multiplier algebra will not coincide with an algebra of bounded Dirichlet series.
Next we will prove that
\begin{equation}\label{eq-4.2}
\m(\h({ces_p}))= \A^{1/q},
\end{equation}
with equality of norms. As explained in the Introduction,  this fact is, to some extent,  analogous
to the case of the space $H(\Disc,ces_p)$, of Taylor series on the unit disc $\Disc$
of the complex plane
with coefficients belonging to $ces_p$, in which case its multiplier algebra  is
the Wiener algebra of  absolutely  convergent Taylor series, which is the smallest algebra
inside $H(\Disc, ces_p)$ containing the  polynomials.


The proof of \eqref{eq-4.2}, which will be
given in Theorem \ref{t-4.8}, is rather
technical. We first discuss the strategy of the proof in order to help its better
understanding. Given $f(s)=\sum_{n=1}^\infty a_nn^{-s}\in\m(\h(ces_p))$, for adequate  values
of  the parameter $\alpha$,
we find a sequence $(g^{m,\alpha})_{m=1}^\infty$ in $\h(ces_p)$ such that
$$
\sum_{n=1}^\infty|a_n|n^{-1/q}=\lim_{{\alpha\to1/q \atop m\to\infty}}
\frac{\Vert fg^{m,\alpha}\Vert_{\h(ces_p)}}{\Vert g^{m,\alpha}\Vert_{\h(ces_p)}}\le\Vert f\Vert_{\m(\h(ces_p))}.
$$
Estimating the norm of $fg^{m,\alpha}$ in $\h(ces_p)$ is complicate since it requires, apart from the Ces\`{a}ro
means, dealing with the coefficients of the product  of two Dirichlet series. With the aim of having
 these coefficients as simple as possible,
we consider functions $g^{m,\alpha}(s)=\sum_{n=1}^\infty b_n^{m,\alpha}n^{-s}$
whose coefficients are supported on certain subsets of the prime numbers  $(p_r)_{r=1}^\infty$.
For an adequate sequence $(r_m)_{m=2}^\infty\subseteq\mathbb{N}$, we require that $b_n^{m,\alpha}\not=0$ only when $n=p_r$ for  $r\ge r_m$.
The key point is that, for coefficients having index of the
form $k=\omega\,p_r$ with $r\ge r_m$ and
$\omega=\prod_{i=1}^{r_m-1}p_i^{t_i}$, $t_1,\dots,t_{r_m-1}\ge0$,
the corresponding coefficient of the product $fg^{m,\alpha}$ is reduced to one term
$$
\sum_{j|k}a_jb_{\frac{k}{j}}^{m,\alpha}=a_\omega b_{p_r}^{m,\alpha}.
$$
In this way, estimating the norms $\Vert fg^{m,\alpha}\Vert_{\h(ces_p)}$
and $\Vert g^{m,\alpha}\Vert_{\h(ces_p)}$ is reduced to estimating  sums of the form
$$
\sum_{{r\ge r_m \atop p_r\le\gamma}}b_{p_r}^{m,\alpha}.
$$
where $\gamma\in[r_m,\infty)$; note that the summation
is taken over the set  $\{r\in\mathbb{N}:\,r\ge r_m \textnormal{ and } p_r\le\gamma\}$.
To this end, we consider the function $\phi(x):=x\log x$ on $[1,\infty)$ and choose, via the Prime Number Theorem,
$(r_m)_{m=2}^\infty\subseteq\N$ such that $p_r$ is sufficiently close to $\phi(r)$ for $r\ge r_m$.
The problem then transformed into  estimating sums of the form
$$
\sum_{r=r_m}^{\phi^{-1}(\gamma)}b_{p_r}^{m,\alpha}.
$$
Finally, good estimates for the above sum are obtained by taking $b_{p_r}^{m,\alpha}=(\phi^\alpha)'(r)$.

We require two lemmata.


\begin{lemma}\label{lemma-J}
Let $0<\beta<1$ and $\phi(x)=x\log x$ for $x\in[1,\infty)$. There exists $x_\beta$ such that for every $r_0\in\mathbb{N}$ with $r_0\ge x_\beta$, $C_1\ge C_2\ge\phi(r_0)$, and $J$ satisfying
\begin{equation}\label{eq-J}
\big\{r\in\mathbb{N}: r\ge r_0 \textnormal{ and }
\phi(r)\le C_2\big\}\subset J\subset\big\{r\in\mathbb{N}: r\ge r_0 \textnormal{ and } \phi(r)\le C_1\big\}
\end{equation}
it follows that
\begin{equation}\label{eq-SumJ}
C_2^\alpha-\phi(r_0)^\alpha
\le
\sum_{r\in J}(\phi^\alpha)'(r)
\le
C_1^\alpha-\phi(r_0-1)^\alpha
\end{equation}
for all $\alpha\le\beta$.
\end{lemma}

\begin{proof}
We consider the Lambert function $W$ on $(0,\infty)$ defined by $W(x)e^{W(x)}=x$;
see \cite{corless-etal}. Then
\begin{equation}\label{eq-4.18}
\phi\Big(\frac{x}{W(x)}\Big)=x.
\end{equation}

Let $r_0\in\mathbb{N}$, $C_1\ge C_2\ge\phi(r_0)$, and $J$ satisfy \eqref{eq-J}.
By \eqref{eq-4.18} and since $\phi$ is increasing and injective  on $[1,\infty)$, we have that $r\le \frac{x}{W(x)}$ if and only if $\phi(r)\le \phi\big(\frac{x}{W(x)}\big)=x$. Then, it follows that
\begin{equation}\label{eq-4.19}
\sum_{r=r_0}^{\big\lfloor\frac{C_2}{W(C_2)}\big\rfloor} h(r)
\le
\sum_{r\in J}h(r)
\le
\sum_{r=r_0}^{\big\lfloor\frac{C_1}{W(C_1)}\big\rfloor} h(r)
\end{equation}
for every positive function $h$. For $\alpha\le\beta$ take $h$ as the derivative  of $\phi^\alpha$, that is,
$$
h(x)=(\phi^\alpha)'(x)=\alpha(x\log x)^{\alpha-1}(\log x+1).
$$
Let $x_\beta$ be sufficiently large so that $h$ is decreasing on $[x_\beta-1,\infty)$.
Such value $x_\beta$ exists as
\begin{eqnarray*}
h'(x) & = & \alpha(\alpha-1)(x\log x)^{\alpha-2}(\log x+1)^2+\alpha(x\log x)^{\alpha-1}\frac{1}{x} \\
& = & \alpha(x\log x)^{\alpha-2}(\log x+1)^2\Big(\alpha-1+\frac{\log x}{(\log x+1)^2}\Big) \\
& \le & \alpha(x\log x)^{\alpha-2}(\log x+1)^2\Big(\,\beta-1+\frac{\log x}{(\log x+1)^2}\Big).
\end{eqnarray*}
Since $\lim_{x\to\infty}\log x(\log x+1)^{-2}=0$ and $\beta-1<0$, there exists $x_\beta$ such that $h'(x)\le 0$
for all $x\ge x_\beta-1$.

Then, for every $M\ge N\ge x_\beta$ it follows that
\begin{align*}
\sum_{r=N}^Mh(r)
\le \sum_{r=N}^M\int_{r-1}^r h(x)\,dx
= \int_{N-1}^M h(x)\,dx
=\phi(M)^\alpha-\phi(N-1)^\alpha
\end{align*}
and
\begin{align*}
\sum_{r=N}^Mh(r)
\ge \sum_{r=N}^M\int_r^{r+1} h(x)\,dx
=\int_N^{M+1} h(x)\,dx
=\phi(M+1)^\alpha-\phi(N)^\alpha.
\end{align*}

Hence, from \eqref{eq-4.19}, if $r_0\ge x_\beta$ we have
\begin{equation*}\label{EQ: hSum2}
\phi\bigg(\bigg\lfloor\frac{C_2}{W(C_2)}\bigg\rfloor+1\bigg)^\alpha-\phi(r_0)^\alpha
\le
\sum_{r\in J}h(r)
\le
\phi\bigg(\bigg\lfloor\frac{C_1}{W(C_1)}\bigg\rfloor\bigg)^\alpha-\phi(r_0-1)^\alpha.
\end{equation*}
From \eqref{eq-4.18},
$$
\phi\bigg(\bigg\lfloor\frac{C_1}{W(C_1)}\bigg\rfloor\bigg)\le\phi\bigg(\frac{C_1}{W(C_1)}\bigg)=C_1
$$
and
$$
\phi\bigg(\bigg\lfloor\frac{C_2}{W(C_2)}\bigg\rfloor+1\bigg)\ge\phi\bigg(\frac{C_2}{W(C_2)}\bigg)=C_2,
$$
and so  \eqref{eq-SumJ} holds.
\end{proof}


\bigskip

Recall that $p_r$ denotes the $r$-th prime number. The Prime Number Theorem
\begin{equation*}\label{eq-PNT}
\lim_{r\to\infty}\frac{p_r}{r\log r}=1,
\end{equation*}
allows to find, for each $2\le m\in\N$,  an integer $r_m\in\N$, with $r_m>m$ and   such that
$$
1-\frac{1}{m}\le\frac{p_r}{r\log r}\le1+\frac{1}{m}, \ \ \textnormal{ for all } r\ge r_m.
$$
Consequently, we have $(r_m)_{m=2}^\infty\subseteq\mathbb{N}$ such that
\begin{equation}\label{eq-4.10}
\frac{mp_r}{m+1}\le r\log r\le\frac{mp_r}{m-1}, \ \ \textnormal{ for all } r\ge r_m.
\end{equation}


\begin{lemma}\label{t-4.9}
Let $\phi(x)=x\log x$ for $x\in[1,\infty)$ and consider the sequence $(r_m)_{m=2}^\infty$ given in \eqref{eq-4.10}. For $1<q<\infty$, there exists $x_q$  such that
\begin{equation*}
\max\Big\{\Big(\frac{m\gamma}{m+1}\Big)^\alpha-\phi(r_m)^\alpha,\,0\Big\}
\le
\sum_{{r\ge r_m\atop p_r\le \gamma}} (\phi^\alpha)'(r)
\le
\Big(\frac{m\gamma}{m-1}\Big)^\alpha-\phi(r_m-1)^\alpha
\end{equation*}
whenever $r_m\ge x_q$, $\gamma\ge p_{r_m}$ and $\alpha\le\frac{1}{q}$.
\end{lemma}

\begin{proof}
Let $\gamma\ge p_{r_m}$ and $\alpha\le\frac{1}{q}$. We apply Lemma \ref{lemma-J} with  $\beta=\frac{1}{q}$, $r_0=r_m$, $C_1=\frac{m\gamma}{m-1}$, $C_2=\frac{m\gamma}{m+1}$ if $\gamma\ge \phi(r_m)\frac{m+1}{m}$, and $C_2=\phi(r_m)$ in other case, and
$$
J=\big\{r\in\mathbb{N}: r\ge r_m \textnormal{ and } p_r\le \gamma\big\}.
$$

We verify that the hypothesis of Lemma \ref{lemma-J} hold.
By \eqref{eq-4.10}  we have  $\phi(r_m)\le\frac{mp_{r_m}}{m-1}\le C_1$ and so $C_1\ge C_2\ge\phi(r_m)$.
The right-hand inclusion in \eqref{eq-J} holds since for every $r\in J$, by \eqref{eq-4.10}, we have that $\phi(r)\le\frac{mp_r}{m-1}\le C_1$.
On the other hand, let $r\ge r_m$ such that $\phi(r)\le C_2$. If $C_2=\frac{m\gamma}{m+1}$, from \eqref{eq-4.10} we have that
$p_r\le\phi(r)\frac{m+1}{m}\le C_2\frac{m+1}{m}=\gamma$ and so $r\in J$. If $C_2=\phi(r_m)$ then $r=r_m\in J$. So, the left-hand inclusion of \eqref{eq-J} holds.

Noting that $C_2^\alpha-\phi(r_m)^\alpha=\max\big\{\big(\frac{m\gamma}{m+1}\big)^\alpha-\phi(r_m)^\alpha,\,0\big\}$, the conclusion follows.
\end{proof}


Now we prove the main result.

\begin{theorem}\label{t-4.8}
For $1<p<\infty$ and $1/p+1/q=1$, we have
$$
\M(\h(ces_p))=\mathcal{A}^{1/q}
$$
with equality of norms.
\end{theorem}

\begin{proof}
Let $f(s)=\sum_{n=1}^\infty a_nn^{-s}\in\m(\h(ces_p))$ and set $a:=(a_n)_{n=1}^\infty$. 
Take $\phi(x)=x\log x$ for $x\in[1,\infty)$, the sequence $(r_m)_{m=2}^\infty$ given in 
\eqref{eq-4.10} and the value $x_q$ provided by Lemma \ref{t-4.9}. For fixed $2\le m\in\N$ 
with $r_m\ge x_q$ and $1/(2q)<\alpha<1/q$, consider the sequence 
$b^{m,\alpha}=(b_n^{m,\alpha})_{n=1}^\infty$ defined by
\begin{equation*}
b_n^{m,\alpha}:=
\left\{\begin{array}{ll}
(\phi^\alpha)'(r) & \textnormal{if $n=p_r$ with  $r\ge r_m$}, \\ 0 & \textnormal{in other case.}
\end{array}\right.
\end{equation*}
Let $g^{m,\alpha}(s):=\sum_{n=1}^\infty b_n^{m,\alpha} n^{-s}$.
Then, from Lemma \ref{t-4.9},
\begin{align}\label{eq-4.16}
\Vert g^{m,\alpha}\Vert_{\h(ces_p)}^p =
\Vert b^{m,\alpha}\Vert_{ces_p}^p
&=\sum_{n=1}^\infty\frac{1}{n^p}\bigg(\sum_{k=1}^n|b_k^{m,\alpha}|\bigg)^p \nonumber
\\ &=
\sum_{n=p_{r_m}}^\infty\frac{1}{n^p}\bigg(\sum_{{r\ge r_m\atop p_r\le n}}(\phi^\alpha)'(r)\bigg)^p \nonumber
\\ &\le
\sum_{n=p_{r_m}}^\infty\frac{1}{n^p}\Big(\frac{mn}{m-1}\Big)^{\alpha p}
\nonumber
\\ &=
\Big(\frac{m}{m-1}\Big)^{\alpha p}\sum_{n=p_{r_m}}^\infty\frac{1}{n^{p(1-\alpha)}}
\nonumber
\\ &\le
\Big(\frac{m}{m-1}\Big)^{\alpha p}\frac{1}{(p(1-\alpha)-1)(p_{r_m}-1)^{p(1-\alpha)-1}}.
\end{align}

We estimate $\Vert fg^{m,\alpha}\Vert_{\h(ces_p)}^p$ from below. Note that for each $k=\omega\,p_r$ 
with $r\ge r_m$ and $\omega=\prod_{i=1}^{r_m-1}p_i^{t_i}$, $t_1,\dots,t_{r_m-1}\ge0$, it follows that
\begin{equation*}
(a\cdot b^{m,\alpha})_k=\sum_{j|k}a_jb_{\frac{k}{j}}^{m,\alpha}=a_\omega b_{p_r}^{m,\alpha}=a_\omega (\phi^\alpha)'(r).
\end{equation*}
Indeed, if $j|k$ and $j\not=\omega$ we have that $k/j\not= p_{\hat{r}}$ for all
$\hat{r}\ge r_m$ and so $b_{k/j}^{m,\alpha}=0$.
Consider the subset of $\mathbb{N}$ given by
$$
\mathcal{P}_m:=
\bigg\{n\in\N:
n=\prod_{i=1}^{r_m-1}p_i^{t_i},\, 0\le t_i\le m
\, \textnormal{ for all } 1\le i\le r_m-1\bigg\}.
$$
Since $\omega p_r=\hat{\omega}p_{\hat{r}}$ with
$\omega,\hat{\omega}\in\mathcal{P}_m$ and $r,\hat{r}\ge r_m$
implies that $\omega=\hat{\omega}$
and $p_r=p_{\hat{r}}$, it follows that the set
$$
\mathcal{O}=\bigcup_{\omega\in\mathcal{P}_m}\,\omega\cdot\big\{p_r:\, r\ge r_m\big\}
$$
is a finite union of disjoint sets.
Then, for any $n_m\ge 3p_{r_m}^{mr_m+1}$, we have that
\begin{eqnarray*}
\Vert fg^{m,\alpha}\Vert_{\h(ces_p)}^p
=
\Vert a\cdot b^{m,\alpha}\Vert_{ces_p}^p
& = &
\sum_{n=1}^\infty\frac{1}{n^p}\bigg(\sum_{k=1}^n|(a\cdot b^{m,\alpha})_k|\bigg)^p
\nonumber
\\ &\ge &
\sum_{n=n_m}^\infty\frac{1}{n^p}
\bigg(\sum_{{k=1 \atop k\in\mathcal{O}}}^n|(a\cdot b^{m,\alpha})_k|\bigg)^p
\nonumber
 \\ & = &
\sum_{n=n_m}^\infty\frac{1}{n^p}
\bigg(\sum_{\omega\in\mathcal{P}_m}
\sum_{{r\ge r_m \atop p_r\le\frac{n}{\omega}}} |(a\cdot b^{m,\alpha})_{\omega\,p_r}|\bigg)^p
 \nonumber
 \\ & = &
\sum_{n=n_m}^\infty\frac{1}{n^p}
\bigg(\sum_{\omega\in\mathcal{P}_m}|a_\omega|
\sum_{{r\ge r_m \atop p_r\le\frac{n}{\omega}}} (\phi^\alpha)'(r)\bigg)^p.
\end{eqnarray*}
Note that $\omega\le p_{r_m}^{mr_m}$ whenever
$\omega\in\mathcal{P}_m$ and so $n/\omega\ge 3p_{r_m}$ for $n\ge n_m$.
Hence, by \eqref{eq-4.10},
$$
\phi(r_m)\le\frac{mp_{r_m}}{m-1}\le\frac{mn}{3(m-1)\omega}\le\frac{mn}{(m+1)\omega}
$$
for every $n\ge n_m$ and $\omega\in\mathcal{P}_m$.
Applying Lemma \ref{t-4.9}, it follows that
$$
\Vert fg^{m,\alpha}\Vert_{\h(ces_p)}^p\ge
\sum_{n=n_m}^\infty\frac{1}{n^p}
\left(\sum_{\omega\in\mathcal{P}_m}|a_\omega|
\bigg(\Big(\frac{mn}{(m+1)\omega}\Big)^\alpha-\phi(r_m)^\alpha\bigg)\right)^p.
$$
Note that if we restrict to $n_m\ge 3p_{r_m}^{mr_m+1+2q}$ we obtain
\begin{eqnarray*}
\Big(\frac{mn}{(m+1)\omega}\Big)^\alpha-\phi(r_m)^\alpha
& = &
\Big(\frac{mn}{(m+1)\omega}\Big)^\alpha
\left(1-\Big(\frac{(m+1)\omega\phi(r_m)}{mn}\Big)^\alpha\right)
\\ & \ge &
\Big(\frac{mn}{(m+1)\omega}\Big)^\alpha\Big(1-\frac{1}{p_{r_m}}\Big)
\end{eqnarray*}
for every $\omega\in\mathcal{P}_m$ and $n\ge n_m$.
Indeed, by \eqref{eq-4.10},
\begin{eqnarray*}
\frac{(m+1)\omega\phi(r_m)}{m}
& \le &
\frac{(m+1)\omega m p_{r_m}}{m(m-1)}
\\ & \le &
\frac{(m+1)p_{r_m}^{mr_m+1}}{m-1}
\\ & \le&
\frac{(m+1)n_m}{3(m-1)p_{r_m}^{2q}}
\le
\frac{\displaystyle  n}{ \displaystyle p_{r_m}^{1/\alpha}},
\end{eqnarray*}
where we use that $\alpha>1/(2q)$. Then,
\begin{eqnarray}\label{eq-4.17}
\Vert fg^{m,\alpha}\Vert_{\h(ces_p)}^p
& \ge &
\sum_{n=n_m}^\infty\frac{1}{n^p}
\left(\sum_{\omega\in\mathcal{P}_m}|a_\omega|\Big(\frac{mn}{(m+1)\omega}\Big)^\alpha
\Big(1-\frac{1}{p_{r_m}}\Big)\right)^p
\nonumber
\\ & = &
\Big(1-\frac{1}{p_{r_m}}\Big)^p
\Big(\frac{m}{m+1}\Big)^{\alpha p}\sum_{n=n_m}^\infty\frac{1}{n^{p(1-\alpha)}}
\bigg(\sum_{\omega\in\mathcal{P}_m}
\frac{|a_\omega|}{\omega^\alpha}\bigg)^p
\nonumber
\\ & \ge &
\Big(1-\frac{1}{p_{r_m}}\Big)^p\Big(\frac{m}{m+1}\Big)^{\alpha p}
\frac{1}{(p(1-\alpha)-1) n_m^{p(1-\alpha)-1}}
\bigg(\sum_{\omega\in\mathcal{P}_m}\frac{|a_\omega|}{\omega^\alpha}\bigg)^p.
\end{eqnarray}
From \eqref{eq-4.16} and \eqref{eq-4.17} it follows that
\begin{align*}
\Vert f\Vert_{\m(\h(ces_p))}^p
&\ge
\frac{\Vert fg^{m,\alpha}\Vert_{\h(ces_p)}^p}{\Vert g^{m,\alpha}\Vert_{\h(ces_p)}^p}
\\ &\ge
\frac{\displaystyle\Big(1-\frac{1}{p_{r_m}}\Big)^p\Big(\frac{m}{m+1}\Big)^{\alpha p}
\frac{1}{(p(1-\alpha)-1) n_m^{p(1-\alpha)-1}}
\bigg(\sum_{\omega\in\mathcal{P}_m}\frac{|a_\omega|}{\omega^\alpha}\bigg)^p}
{\displaystyle \Big(\frac{m}{m-1}\Big)^{\alpha p}\frac{1}{(p(1-\alpha)-1)(p_{r_m}-1)^{p(1-\alpha)-1}}}
\\ &=
\frac{(p_{r_m}-1)^{p(2-\alpha)-1}}{p_{r_m}^p}\Big(\frac{m-1}{m+1}\Big)^{\alpha p}
\frac{1}{n_m^{p(1-\alpha)-1}}
\bigg(\sum_{\omega\in\mathcal{P}_m}\frac{|a_\omega|}{\omega^\alpha}\bigg)^p.
\end{align*}
Taking limit as $\alpha\to1/q$ we have
$$
\Vert f\Vert_{\m(\h(ces_p))}^p\ge
\frac{(p_{r_m}-1)^p}{p_{r_m}^p}\Big(\frac{m-1}{m+1}\Big)^{{p}/{q}}
\bigg(\sum_{\omega\in\mathcal{P}_m}\frac{|a_\omega|}{\omega^{1/q}}\bigg)^p.
$$
Finally, making $m\to\infty$ we conclude
$$
\Vert f\Vert_{\m(\h(ces_p))}^p\ge
\bigg(\sum_{\omega\in\N}\frac{|a_\omega|}{\omega^{1/q}}\bigg)^p.
$$
\end{proof}


\section{Further facts on  multipliers on $\h(ces_p)$}\label{section-4}


First we study the compactness of the multipliers on $\h(ces_p)$. It turns out that
there is no other compact multiplier than zero.

\begin{theorem}\label{t-5.1}
Let $f\in\m(\h({ces_p}))$. Suppose that the associated operator
$$
g\in\h(ces_p)\mapsto M_f(g):=f g\in\h(ces_p)
$$
is compact. Then  $f=0$.
\end{theorem}

\begin{proof}
Consider the sequence $\{m^{1/q}m^{-s}\}_{m=1}^\infty$ in $\h(ces_p)$.
It is bounded  as, for $m\ge2$, we have that
$$
\Vert m^{1/q}m^{-s}\Vert_{\h(ces_p)}  =  m^{1/q}\Vert e^m\Vert_{ces_p}
=m^{1/q}\Big(\sum_{n=m}^\infty\frac1{n^p}\Big)^{1/p}
\le
\frac{2^{1/q}}{(p-1)^{1/p}}.
$$
Then, by compactness of $M_f$,
there exists a subsequence $\{m_k^{1/q}m_k^{-s}\}_{k=1}^\infty$
such that $\{M_f(m_k^{1/q}m_k^{-s})\}_{k=1}^\infty$ converges in norm to some
$g\in\h({ces_p})$.
For $s_0\in\mathbb{C}_{1/q}$, since the point evaluation $\delta_{s_0}$ is bounded on
$\h(ces_p)$, we have
$$
\delta_{s_0}\Big(M_f(m_k^{1/q}m_k^{-s})\Big)\xrightarrow[k\rightarrow\infty]{}
\delta_{s_0}(g)=g(s_0).
$$
On the other hand,
$$
\delta_{s_0}\Big(M_f(m_k^{1/q}m_k^{-s})\Big)
=f(s_0)m_k^{1/q-s_0}\xrightarrow[k\rightarrow\infty]{} 0.
$$
Thus, $g=0$. Hence, $\{M_f(m_k^{1/q}m_k^{-s})\}_{k=1}^\infty$
converges to zero in the norm of $\h(ces_p)$.

We estimate from below
$\Vert M_f(m_k^{1/q}m_k^{-s})\Vert_{\h(ces_p)}
=m_k^{1/q}\Vert M_f(m_k^{-s})\Vert_{\h(ces_p)}$.
Let $f(s)=\sum_{n=1}^\infty a_nn^{-s}$.
We have seen in the proof of Theorem \ref{t-4.5} that
$$
\Vert M_f(m^{-s})\Vert_{\h(ces_p)}^p
=\Vert m^{-s}f\Vert_{\h(ces_p)}^p
=\sum_{j=1}^\infty\Big(\sum_{i=1}^j|a_i|\Big)^p\sum_{n=jm}^{(j+1)m-1}\frac{1}{n^p}.
$$
Since
$$
\sum_{n=jm}^{(j+1)m-1}\frac{1}{n^p}\ge\frac{m}{\big((j+1)m-1\big)^p}
\ge\frac{m}{(2jm)^p},
$$
it follows that
$$
\Vert M_f(m^{-s})\Vert_{\h(ces_p)}^p
\ge\frac{m}{(2m)^p}\sum_{j=1}^\infty\frac{1}{j^p}\Big(\sum_{i=1}^j|a_i|\Big)^p=
\frac{m}{(2m)^p}\Vert f\Vert_{\h(ces_p)}^p.
$$
Then,
$$
\Vert M_f(m_k^{1/q}m_k^{-s})\Vert_{\h(ces_p)}
\ge m_k^{1/q}\frac{m_k^{1/p}}{2m_k}\Vert f\Vert_{\h(ces_p)}=
\frac{1}{2}\Vert f\Vert_{\h(ces_p)}.
$$
Taking $k\to\infty$ we have that $\Vert f\Vert_{\h(ces_p)}\le0$ and so, $f=0$.
\end{proof}


Next we discuss how ``close'' is the space $\h(ces_p)$ to its multiplier algebra.
Let us first note that $\m(\h(ces_p))=\A^{1/q}\varsubsetneq\h(ces_p)$.
Indeed, in other case   the point evaluation $\delta_{1/q}$ at the point $s_0=1/q$, which
belongs to the dual space of $\A^{1/q}$, would belong to $\h(ces_p)^*$.
But this is not true as $\delta_{1/q}\in\h(ces_p)^*$ is precisely  $(n^{-1/q})_{n=1}^\infty\in d(q)$ (see the proof of Theorem \ref{t-3.4}), which corresponds to
$(n^{-1/q})_{n=1}^\infty\in \ell^q$.

The multiplier algebra $\m(\h(ces_p))$ is ``close'' to $\h(ces_p)$
in the sense shown by the following example.
For $f(s)=\sum_{n=1}^\infty a_n n^{-s}\in\h(ces_p)$ and $\varepsilon>0$, set
$$
f_{\varepsilon}(s):=\sum_{n=1}^\infty\frac{a_nn^{-\varepsilon}}{n^s}.
$$
 Theorem \ref{t-3.3} shows that  $\sigma_a(\h({ces_p}))=1/q$. Then
$$
\sum_{n=1}^\infty\frac{|a_nn^{-\varepsilon}|}{n^{1/q}}
=\sum_{n=1}^\infty\frac{|a_n|}{n^{1/q+\epsilon}}<\infty .
$$
Consequently, $f_{\varepsilon}\in\A^{1/q}=\m(\h({ces_p}))$.
The question arises: for which sequences $(b_n)_{n=1}^\infty$
it is the case that $\sum_{n=1}^\infty a_nb_nn^{-s}\in\m(\h(ces_p))$ whenever
$\sum_{n=1}^\infty a_n n^{-s}\in \h(ces_p)$?
Recall that these sequences are called the Schur multipliers between
$\h(ces_p)$  and $\A^{1/q}$.

\begin{theorem}\label{t-5.2}
A sequence $(b_n)_{n=1}^\infty$ satisfies that for every $\sum_{n=1}^\infty a_n n^{-s}\in \h(ces_p)$
the series $\sum_{n=1}^\infty a_nb_nn^{-s}$ is a multiplier on $ \h(ces_p)$ if and only if
$$
\big(b_n n^{-1/q}\big)_{n=1}^\infty \in d(q),
$$
where $d(q)$ is the dual space to $ces_p$, that is, the following condition holds
$$
\sum_{n=1}^\infty \sup_{k\ge n}\bigg(\frac{|b_k|^q}{k}\bigg)<\infty.
$$
\end{theorem}

\begin{proof}
Denote by $\ell_q^1$ the Banach space of complex sequences $a=(a_n)_{n=1}^\infty$ such that $\Vert a\Vert_{\ell_q^1}:=\sum_{n=1}^\infty|a_n|n^{-1/q}<\infty$.
The sequence $b=(b_n)_{n=1}^\infty$ being a Schur  multiplier between
$\h(ces_p)$ and $\A^{1/q}$ corresponds precisely to the operator $T_b$, defined by
$$
a=(a_n)_{n=1}^\infty\in ces_p\mapsto
T_b(a):=(a_nb_n)_{n=1}^\infty\in\ell_q^1,
$$
being well defined and, via the closed graph theorem, bounded.

Suppose that $T_b$ is well defined. For every $a=(a_n)_{n=1}^\infty\in ces_p$ we have that
$$
\Big|\Big\langle\big(b_n n^{-1/q}\big)_{n=1}^\infty,(a_n)_{n=1}^\infty\Big\rangle\Big|\le\sum_{n=1}^\infty\frac{|a_nb_n|}{n^{1/q}}=\Vert T_b(a)\Vert_{\ell_q^1}\le\Vert T_b\Vert\cdot\Vert a\Vert_{ces_p}
$$
and so $\big(b_n n^{-1/q}\big)_{n=1}^\infty\in ces_p^*=d(q)$.

Conversely, suppose that $\big(b_n n^{-1/q}\big)_{n=1}^\infty\in ces_p^*=d(q)$. For every $a=(a_n)_{n=1}^\infty\in ces_p$ we have that
$$
\bigg|\sum_{n=1}^\infty\frac{a_nb_n}{n^{1/q}}\bigg|=\Big|\Big\langle\big(b_n n^{-1/q}\big)_{n=1}^\infty,(a_n)_{n=1}^\infty\Big\rangle\Big|\le\Big\Vert \big(b_n n^{-1/q}\big)_{n=1}^\infty\Big\Vert_{ces_p^*}\cdot\Vert a\Vert_{ces_p}.
$$
Let $c=(c_n)_{n=1}^\infty$ with
$c_n:=\frac{|a_n|\overline{b_n}}{|b_n|}$ if $b_n\not=0$
and $c_n=0$ in other case. Then
$|c_n|\le|a_n|$ for all $n\ge1$. Since $a\in ces_p$,
it follows that $c\in ces_p$ and $\|c\|_{ces_p}\le\|a\|_{ces_p}$.
Hence,
\begin{align*}
\sum_{n=1}^\infty\frac{|a_nb_n|}{n^{1/q}}=\bigg|\sum_{n=1}^\infty\frac{c_nb_n}{n^{1/q}}\bigg|&\le\Big\Vert \big(b_n n^{-1/q}\big)_{n=1}^\infty\Big\Vert_{ces_p^*}\cdot\Vert c\Vert_{ces_p}
\\ & \le
\Big\Vert\big(b_n n^{-1/q}\big)_{n=1}^\infty\Big\Vert_{ces_p^*}\cdot\Vert a\Vert_{ces_p}.
\end{align*}
So $T_b(a)\in\ell_q^1$.
\end{proof}


\begin{example}
For $\alpha>\frac{1}{q}$ set $b_n=(\log n)^{-\alpha}$, for $n\ge2$. Then
$$
\sum_{n=2}^\infty \sup_{k\ge n}\bigg(\frac{|(\log k)^{-\alpha}|^q}{k}\bigg)
=\sum_{n=2}^\infty \frac{1}{n(\log  n)^{q\alpha}}<\infty.
$$
Thus,  for every  $\sum_{n=1}^\infty a_n n^{-s}\in \h(ces_p)$ we have that
$$
\sum_{n=2}^\infty\frac{a_n}{n^s\log^\alpha n}
$$
is a multiplier on $\h(ces_p)$.
\end{example}



\end{document}